\documentclass[fleqn]{article}
\usepackage{amssymb,latexsym, graphicx, xcolor}

\def\<{\langle}
\def\>{\rangle}

\newcommand{\dm}{\mathsf{dim}}
\newtheorem{res}{Result}[section]
 
\newtheorem{lemma}[res]{Lemma} 
\newtheorem{cor}[res]{Corollary} 
\newtheorem{prop}[res]{Proposition} 
\newtheorem{theo}[res]{Theorem} 
\newtheorem{prob}[res]{Problem} 
\newtheorem{quest}[res]{Question} 
 
\newtheorem{conj}[res]{Conjecture} 
\newtheorem{note}{Remark}

\title{Computations regarding certain graphs associated to finite polar spaces}
\author{Antonio Pasini}
\begin{document}
\maketitle

\begin{abstract}
We consider various regular graphs defined on the set of elements of given rank of a finite polar space. It is likely that no two such graphs, of the same kind but defined for different ranks, can have the same degree. We shall prove this conjecture under the hypothesis that the considered rank are not too small.   
\end{abstract} 

\section{Introduction}

We presume that the reader is acquainted with basics on polar spaces; if not, we refer him$/$her to \cite[Chapter 7]{Tits}, also \cite[Chapter 7]{BC}. All polar spaces to be considered here are non-degenerate and thick-lined, namely all of their lines admit at least three points. We recall that the singular subspaces of a (non-degenerate) polar space are projective spaces. 

Following a well establised custom, we denote collinearity by the symbol $\perp$. Thus, given two sets of points $X$ and $Y$, if all points of $X$ are collinear with all points of $Y$ then we write $X\perp Y$ and $X^\perp$ is the set of points collinear with all points of $X$. For two singular subspaces $X$ and $Y$, we put $X^Y := \langle X, Y\cap X^\perp\rangle$. Clearly, $X^Y$ is the smallest one among the singular subspaces which contain $X$ and admit the largest intersection with $Y$. Note that $\dm X^Y = \dm Y^X$.  

\subsection{A motivation for this paper}

Let $\Delta$ be a polar space of rank $n \geq 3$ and for $i = 0, 1,..., n-1$ let $\Delta_i$ be the set of $i$-dimensional singular subspaces of $\Delta$. For $0 \leq i < j \leq \ell \leq n-1$ and $-1 \leq k \leq i$ such that $i+j-k = \ell$, let $\Gamma_{i,j;k,\ell}(\Delta)$ be the bipartite graph with $\Delta_i\cup\Delta_j$ as the set of vertices and the edges of which are the pairs $\{I,J\}$ with $I\in \Delta_i$, $J\in \Delta_j$ and such that $\dm(I\cap J) = k$ and $\dm I^J = \dm J^I = \ell$.  

In \cite{DSVM} De Schepper and Van Maldeghem consider the following problem:

\begin{prob}\label{prob 1}
Can we recover $\Delta$ from $\Gamma_{i,j;k,\ell}(\Delta)$? 
\end{prob} 
\begin{note}
{\rm The class of graphs considered in \cite{DSVM} is slightly more general than defined here. Indeed the definition of the adjacency relation of the graph $\Gamma_{i,j;k,\ell}(\Delta)$ as stated above can be re-phrased in terms of Weyl distances (see \cite{DSVM}) and, in that setting, it also makes sense for graphs which are not covered by our formulation. However, our definition is sufficient to give the reader a right idea of the problem studied in \cite{DSVM}.}
\end{note} 

De Schepper and Van Maldeghem \cite{DSVM} answer Problem \ref{prob 1} in the affirmative but for the case where $i < \ell = n-1$, which cannot be treated by the techinque they exploit in \cite{DSVM}. A subcase of this wild case is studied by De Schepper and the author in \cite{DSP}. Explicitly, assuming that $i+j -k = \ell = n-1$, an affirmative answer is obtained for Problem \ref{prob 1}. In particular it is proved in \cite{DSP} that, when $i+j - k = n-1 = \ell$, Problem \ref{prob 1} can be answered in the affirmative provided that the following is true:  

\begin{theo}\label{theo 1}
Let $i+j-k = n-1$ with $i < j < n-1$ Then no automorphism of the bipartite graph $\Gamma_{i,j;k,n-1}(\Delta)$ permutes the two classes $\Delta_i$ and $\Delta_j$ of its bipartition. 
\end{theo} 

Theorem \ref{theo 1} is proved in \cite{DSP} by purely geometrical arguments, valid when the lines of $\Delta$ have at least four points. 
When the lines of $\Delta$ have size 3 a different approach is necessary, relying on computations. 

In this paper we shall perform those computations, in a more general setting than \cite{DSP} would require. Indeed we only assume finiteness while in view of \cite{DSP} we could restrict ourselves to the case where the lines of $\Delta$ have just three points. However, treating the general case, with lines of size $s+1$ instead of $3$, makes no difference: it just amounts to write $s$ instead of $2$. 

It goes without saying that the proof of Theorem \ref{theo 1} in \cite{DSP} rests upon some of the results to be proved in the sequel. Thus, this paper can be regarded as a completion of \cite{DSP}; even better, it stands in support of \cite{DSP}. 

\subsection{A few graphs defined on $\Delta_i$}       

Henceforth $\Delta$ is a finite polar space of rank $n > 2$ and order $(s, s,..., s, t)$, where $s$ and $t$ are positive integers, $s > 1$ and $t \in \{1, s^{1/2}, s, s^{3/2}, s^2\}$. 

A number of `natural' symmetric relations can be considered on the set $\Delta_i$, which define a graph on it. The collinearity relation $\perp$ is one of them. Here is another one: for $X, Y\in \Delta_i$, we write $X\approx Y$ if $\dm(X\cap Y) = i-1$. We can also consider the join $\perp\cup\approx$ and the meet of $\perp\cap\approx$ of $\perp$ and $\approx$. The latter is the collinarity relation of the $i$-grassmannian of $\Delta$ for $i < n-1$ while $\approx$ is the collinearity relation in the $(n-1)$-grassmanian of $\Delta$. Note also that $(\Delta_0, \approx)$ is a complete graph and $(\Delta_{n-1}, \perp)$ is totally disconnected (no edges at all). 

One more graph is considered in \cite{DSP}, namely $(\Delta_i, \perp_{\mathrm{max}})$ where $X \perp_{\mathrm{max}} Y$ if $X\perp Y$ and $\dm\langle X, Y\rangle = n-1$. Clearly, if $i < (n-2)/2$ then $(\Delta_i, \perp_{\mathrm{max}})$ is totally disconnected.  Also, $\perp$ and $\perp_{\mathrm{max}}$ coincide on $\Delta_{n-2}$.

For the rest of this subsection we assume that $i < j < n-1$, as in Theorem \ref{theo 1}. For $u\in \{i,j\}$ and $X, Y\in \Delta_u$, put $\mathrm{N}(X,Y) := \{X,Y\}^{\sim}$, where $\sim$ stands for the adjacency relation of $\Gamma_{i,j;k,n-1}$. 

We recall that, for two distinct non-collinear points $x, y$ of $\Delta$, the set $\{x,y\}^{\perp\perp}$ is called a hyperbolic line. All hyperbolic lines of $\Delta$ have the same size, say $h(\Delta)$. Clearly, $h(\Delta) \geq 2$, since in any case $\{x,y\}^{\perp\perp}\supseteq \{x,y\}$. We call $h(\Delta)$ the {\em hyperbolic order} of $\Delta$. The following is proved in \cite{DSP}: 

\begin{prop}\label{N 1}
Let $X,Y\in \Delta_u$, $u\in \{i,j\}$ and $X\neq Y$. 
\begin{itemize}
\item[$(1)$] If $h(\Delta) > 2$ then $\mathrm{N}(\mathrm{N}(X,Y)) = \{X,Y\}$ if and only if $X\perp Y$.
\item[$(2)$] Let $h(\Delta) = 2$. Then $\mathrm{N}(\mathrm{N}(X,Y)) = \{X,Y\}$ if and only if either $X\perp Y$ or $X\approx Y$.
\end{itemize}
\end{prop}

\begin{cor}\label{N 1 bis}
If an automorphism of $\Gamma_{i,j;k,n-1}$ switches $\Delta_i$ and $\Delta_j$ then it induces an isomorphism between either $(\Delta_i, \perp)$ and $(\Delta_j,\perp)$ (if $h(\Delta) > 2$) or $(\Delta_i, \perp\cup\approx)$ and $(\Delta_j,\perp\cup\approx)$ (if $h(\Delta) = 2$). 
\end{cor} 

With $u\in \{i,j\}$ and $X, Y \in \Delta_u$, $X\neq Y$, we write $X\sim_u Y$ if $\mathrm{N}(\mathrm{N}(X,Y)) = \{X,Y\}$. The following is also proved in \cite{DSP}.

\begin{prop}\label{N 2}
With $\{u,v\} = \{i,j\}$, let $X, Y \in \Delta_u$, $X\neq Y$. 
\begin{itemize}
\item[$(1)$] Let $u \geq (n-2)/2$. If $u = n-2$, assume moreover that $h(\Delta) > 2$. Then $X\perp_{\mathrm{max}} Y$ if and only if $X\sim_u Y$ and $\mathrm{N}(X,Y)$ is a non-trivial clique of $(\Delta_v,\sim_v)$.
\item[$(2)$] Let $u < (n-2)/2$. Then $X\sim_u Y$ only if the set $\mathrm{N}(X,Y) \subseteq \Delta_v$ is not a clique of the graph $(\Delta_v,\sim_v)$. 
\item[$(3)$] Let $u = n-2$ and $v = 0$. Suppose also that $h(\Delta) = 2$. Then $X\approx Y$ if and only if $X\sim_u Y$ and $\mathrm{N}(X,Y)$ is a non-trivial clique of $(\Delta_v,\sim_v)$. 
\item[$(4)$] Let $u = n-2$ (hence $u = j > i = v$) but $v > 0$. Assume moreover that $h(\Delta) = 2$ and suppose that $X\sim_u Y$ but $X\not\perp Y$. Then $\mathrm{N}(X,Y)$ is not a clique in $(\Delta_v,\sim_v)$. 
\end{itemize}
\end{prop}
Turning back to Theorem \ref{theo 1}, note that the conditions $i < j$ and $i+j-k = n-1$ force $j > (n-2)/2$. By this fact and Proposition \ref{N 2} we immediately get the following corollaries:

\begin{cor}\label{N 2 bis}
If $i < (n-2)/2$ then no automorphism of $\Gamma_{i,j;k,n-1}(\Delta)$ can switch $\Delta_i$ and $\Delta_j$.
\end{cor}

\begin{cor}\label{N 2 ter}
Let $(n-2)/2 \leq i$. If an automorphism of $\Gamma_{i,j;k,n-1}(\Delta)$ switches $\Delta_i$ and $\Delta_j$, it induces an isomorphism from $(\Delta_i, \perp_{\mathrm{max}})$ to $(\Delta_j, \perp_{\mathrm{max}})$.
\end{cor}

In view of Corollaries \ref{N 1 bis}, \ref{N 2 bis} and \ref{N 2 ter}, the following is sufficient for a proof of Theorem \ref{theo 1} in the finite case:

\begin{theo}\label{theo 2}
Let $(n-2)/2 \leq i < j \leq n-1$ and let $\Delta$ be finite. Then $(\Delta_i, \perp)\not\cong(\Delta_j, \perp)$, $(\Delta_i, \perp\cup\approx)\not\cong(\Delta_j, \perp\cup\approx)$ and $(\Delta_i, \perp_{\mathrm{max}})\not\cong(\Delta_j, \perp_{\mathrm{max}})$.
\end{theo}

We shall obtain this theorem in Section \ref{proof} as a trivial consequence of results on the degrees of the graphs considered in it, to be proved in Section \ref{degree sec}. 

\section{Comparing $|\Delta_i|$ and $|\Delta_j|$}

It goes without saying that if we could answer the following question in the negative then Theorem \ref{theo 1} would be proved:

\begin{quest}\label{quest}
Can we have $|\Delta_i| = |\Delta_j|$ for suitable choices of the types $i < j$ and the orders $s$ and $t$? 
\end{quest}

Regretfully, Question \ref{quest} admits an affirmative answer when $s = t$ with $2n\equiv 1~(\mathrm{mod}~3)$. Explicitly, with $s$ and $t$ chosen in this way, if $i = (2n-4)/3$ and $j = (2n-1)/3$ then $|\Delta_i| = |\Delta_j|$ (see the second part of Proposition \ref{maximum}, next subsection). This is the unique example we are aware of, but the outcome of the following elementary investigation of the mapping $i \rightarrow |\Delta_i|$ warns us that the equality $|\Delta_i| = |\Delta_j|$ could possibly hold in many other cases.   

\subsection{The function $\delta_{n,s,t}:i\rightarrow |\Delta_i|$}\label{max sec}

Let $\delta_{n,s,t}$ be the function which maps $i \in \{0,1,..., n-1\}$ onto the cardinality $\delta_{n,s,t}(i) = |\Delta_i|$ of $\Delta_i$. By a standard double counting we get
\begin{equation}\label{2 old}
\frac{\delta_{n,s,t}(i+1)}{\delta_{n,s,t}(i)} ~ = ~ \frac{(s^{n-(i+1)-1}t+1)(s^{n-(i+1)}-1)}{(s^{(i+1)+1}-1)},
\end{equation}
for $i = 0, 1,..., n-2$. Equation (\ref{2 old}) combined with the well known equality $\delta_{n,s,t}(0) = (s^{n-1}t+1)(s^n-1)/(s-1)$ implies
\begin{equation}\label{1 old}
\delta_{n,s,t}(i) ~ = ~ \frac{\prod_{u=0}^i(s^{n-u-1}t+1)\cdot\prod_{u=0}^i(s^{n-u}-1)}{\prod_{u=0}^i(s^{u+1}-1)}.
\end{equation}
Define $i_{\mathrm{max}}$ as follows:
\[i_{\mathrm{max}} ~ := ~ \left\{\begin{array}{ll}
 \lfloor (2n-2+\log_st)/3\rfloor & \mbox{if } t \leq s,\\
 \lceil (2n-2+\log_st)/3\rceil - 1 & \mbox{if } s < t.
\end{array}\right.\]

\begin{prop}\label{maximum}
If either $s \neq t$ or $s = t$ but $3$ does not divide $2n-1$, then the function $\delta_{n,s,t}$ attains its maximum value at $i_{\mathrm{max}}$, it is increasing in $\{0, 1,..., i_{\mathrm{max}}\}$ and decreasing in $\{ i_{\mathrm{max}},..., n-1\}$. 

On the other hand, let $s = t$ and $2n\equiv 1~(\mathrm{mod}~3)$. Then $\delta_{n,s,s}(i_{\mathrm{max}}-1) = \delta_{n,s,s}(i_{\mathrm{max}})$. In this case $\delta_{n,s,s}$ attains its maximum at $i_{\mathrm{max}}$ as well as $i_{\mathrm{max}}-1$, it is increasing in $\{0, 1,..., i_{\mathrm{max}}-1\}$ and decreasing in $\{ i_{\mathrm{max}},..., n-1\}$. 
\end{prop}
{\bf Proof.} Suppose $s\neq t$ or $2n\not\equiv 1~(\mathrm{mod}~3)$. It is straightforward to check that the ratio $\delta_{n,s,t}(i+1)/\delta_{n,s,t}(i)$ is greater (smaller) than $1$ if and only if $i$ is smaller (respectively, greater or equal) than $i_{\mathrm{max}}$. When $s = t$ and $2n\equiv 1~(\mathrm{mod}~3)$, then $\delta_{n,s,s}(i+1)/\delta_{n,s,s}(i)$ is greater (smaller) than $1$ if and only if $i$ is smaller (greater) than $i_{\mathrm{max}}-1$. The statement follows.  \hfill $\Box$\\

In short, $\delta_{n,s,t}$ is increasing for approximately the first two thirds of its domain and decreasing in the last third.  

\begin{note}
{\rm If $(2n-2+\log_s t)/3$ is non integral then $\lfloor (2n-2+\log_st)/3\rfloor = \lceil (2n-2+\log_st)/3\rceil -1$. In this case the two parts of the definition of $i_{\mathrm{max}}$ could be fused in one single formula. However, when $\log_st$ is integral and $3$ divides $(2n-2+\log_st)$ then $\lfloor (2n-2+\log_st)/3\rfloor = (2n-2+\log_st)/2 > (2n-2+\log_st)/3-1 = \lceil (2n-2+\log_st)/3\rceil -1$. In this case me must stick to the above definition as stated.  If we like, we can replace the conditions $t\leq s$ and $s < t$ of that definition with $t < s$ and $s\leq t$ respectively, but we should accordingly modify the second part of Proposition \ref{maximum}.}
\end{note} 

\subsection{A few results in the general case}

In this subsection we shall obtain two (admittedly lax) necessary conditions to be satisfied by $i$ and $j$ for the equality $|\Delta_i| = |\Delta_j|$ to hold. 

Let $j > i$, as in Question \ref{quest}. By (\ref{1 old}), we have $|\Delta_i| \geq |\Delta_j|$ if and only if 
\begin{equation}\label{1 new}
\frac{\prod_{u= i+1}^j(s^{n-u-1}t+1)(s^{n-u}-1)}{\prod_{u = i+1}^j(s^{u+1}-1)} ~ \leq ~ 1,
\end{equation}
namely $\prod_{u= i+1}^j(s^{n-u-1}t+1)(s^{n-u}-1) \leq  \prod_{u = i+1}^j(s^{u+1}-1)$. Equivalently:
\[\prod_{u= i+1}^j((s^{n-u-1}t+1)\cdot\sum_{v=0}^{n-u-1}s^v) ~ \leq ~ \prod_{u = i+1}^j\sum_{v=0}^{u}s^v.\]
Suppose that this inequality holds. Then $\prod_{u = i+1}^js^{2n-2u-2+\log_st}  <  \prod_{u= i+1}^js^{u+1}$, 
namely $\sum_{u = i+1}^j(2n-2u-2+\log_st)  <  \sum_{u=i+1}^j(u+1)$. 
Equivalently,
\[\sum_{u = i+1}^j(2n+\log_st) ~ < ~ 3\cdot\sum_{u=i+1}^j(u+1).\]
Therefore $(j-i)(2n +\log_st)  < 3\cdot((j-i)(i+1) + ((j-i)(j-i+1))/2)$, namely $2n +\log_st < 3(i+j+3)/2$. Thus, we have proved the following: 

\begin{prop}\label{lemma com 1}
If $|\Delta_i| \geq |\Delta_j|$ then $2n +\log_st < 3(i+j+3)/2$. 
\end{prop} 

Suppose now that $|\Delta_i| \leq |\Delta_j|$, namely
\begin{equation}\label{2 new}
\frac{\prod_{u= i+1}^j(s^{n-u-1}t+1)(s^{n-u}-1)}{\prod_{u = i+1}^j(s^{u+1}-1)} ~ \geq ~ 1.
\end{equation}
Then $\prod_{u= i+1}^j((s^{n-u-1}t+1)\cdot\sum_{v=0}^{n-u-1}s^v) \geq \prod_{u = i+1}^j\sum_{v=0}^{u}s^v$. This inequality implies that 
$\prod_{u = i+1}^js^{2n-2u+\log_st} >  \prod_{u= i+1}^js^{u}$. Equivalently,
\[\sum_{u = i+1}^j(2n-2u+\log_st) ~ > ~ \sum_{u=i+1}^ju.\]
By easy manipulations now we obtain: 

\begin{prop}\label{lemma com 2}
If $|\Delta_i| \leq |\Delta_j|$ then $2n +\log_st > 3(i+j+1)/2$. 
\end{prop} 

By combining Propositions \ref{lemma com 1} and \ref{lemma com 2} and recalling that $2\log_st$ is integral, we easily obtain the following:

\begin{cor}\label{res com 1}
If $|\Delta_i| = |\Delta_j|$ then
\begin{equation}\label{com 3}
2n + \log_st - 2 ~\geq ~ \frac{3}{2}(i+j) ~ \geq ~ 2n + \log_st - 4.
\end{equation}
\end{cor} 
In simple words, if $|\Delta_i| = |\Delta_j|$ then the average $(i+j)/2$ of $i$ and $j$ is close to $i_{\mathrm{max}}\approx 2n/3$.  

\begin{prop}\label{lemma com 3}
If $|\Delta_i| \leq |\Delta_j|$ then $2n-4+\log_st \geq 2i+j$. 
\end{prop}
{\bf Proof.} Let $|\Delta_i|\leq |\Delta_j|$. By way of contradiction, suppose that $2n+\log_st < 2i+j + 4$. However $2n +\log_st \geq 3(i+j)/2 + 1$ by Proposition \ref{lemma com 2} (recall that $2\log_st$ is integral). It follows that $2i+j+3 > 3(i+j)/2$, namely $4i + 2j + 6 > 3i + 3j$. Equivalently, $i \geq j-5$. However, $i < j$ by assumption. Therefore $j-5 \leq i \leq j-1$. However
\begin{equation}\label{1 new prel bis}
\prod_{u= i+1}^j(s^{n-u-1}t+1)(s^{n-u}-1) ~ \geq ~ \prod_{u = i+1}^j(s^{u+1}-1).
\end{equation} 
since $|\Delta_i| \leq |\Delta_j|$ by assumption. It is straightforward to check that inequality (\ref{1 new prel bis}) fails to hold if $j-4 \leq i \leq j-1$ and $2n+\log_st < 2i+j + 4$. We have reached a contradiction, as we wanted. 

The above proof rests on Proposition \ref{lemma com 2}, but we can also give a proof that does not depend on that proposition. Here it is. Rewrite (\ref{1 new prel bis}) as follows:
\begin{equation}\label{1 new prel ter}
\prod_{v= n-j}^{n-i-1}(s^{v-1}t+1)(s^v-1) ~ \geq ~ \prod_{u = i+1}^j(s^{u+1}-1).
\end{equation}
The factors of the product at the left side of (\ref{1 new prel ter}) increase more quickly than those at the right side. If $j > i+1$ then (\ref{1 new prel ter}) holds only if the largest factor of the left product is greater than the largest factor at right. Explicitly,
\begin{equation}\label{1 - A}
(s^{n-i-2}t+1)(s^{n-i-1}-1) ~ \geq ~ s^{j+1}-1
\end{equation}
with strict inequality when $j > i+1$. By standard manipulation, from (\ref{1 - A}) we obtain the inequality $2n-4+\log_st \geq 2i+j$.  \hfill $\Box$ \\

In similar ways we can also prove the following:

\begin{prop}\label{lemma com 4}
If $|\Delta_i| \leq |\Delta_j|$ then $2n-3+\log_st \leq i+2j$. 
\end{prop}

The next corollary immediately follows from Propositions \ref{lemma com 3} and \ref{lemma com 4}. 

\begin{cor}
If $|\Delta_i| = |\Delta_j|$ then $2i+j+4 \leq 2n+\log_st \leq i+2j+3$.
\end{cor} 
 
\subsection{A special case}

A particular attention is payed in \cite{DSP} to the case $(n-2)/2\leq i < j = n-2$. This section is devoted to that case.  

Put $j = n-2$ and assume that $|\Delta_i| = |\Delta_{n-2}|$. Then (\ref{com 3}) yields 
\begin{equation}\label{com spec 1}
n + 2\log_st + 2 ~ \geq ~ 3i ~ \geq ~ n + 2\log_st - 2.
\end{equation} 

\begin{prop}\label{res com 2}
Let $(n-2)/2 \leq i < n-2$. Then $|\Delta_i| = |\Delta_{n-2}|$ if and only if $n = 5$, $i = 2$ and $t = s$.
\end{prop} 
{\bf Proof.} We have $\lceil n/2\rceil -1 \leq i \leq \lfloor n/2\rfloor$ since $(n-2)/2\leq i$ (by assumption) and $i\leq n/2$ (by Proposition \ref{lemma com 3} with $j = n-2$). In other words, if $n$ is even, say $n = 2m$, then $i \in \{m-1, m\}$. If $n$ is odd, say $n = 2m+1$, then $i = m$. We shall split the proof in two parts, considering the case $n = 2m$ first; next we turn to the case $n = 2m+1$. 

\medskip

\noindent
{\bf Even case: $n = 2m$.} Suppose firstly that $i = m-1$. The first inequality of (\ref{com spec 1}) (which follows from the hypothesis that $|\Delta_i| \leq |\Delta_j|$) yields $m  \leq  2\log_s t + 5$ ($\leq 9$). In order to avoid to examine too many cases, we shall firstly prove that, if $m \geq 6$ then $|\Delta_{m-1}| > |\Delta_{n-2}|$, namely
\begin{equation}\label{1 bis}
\prod_{u=m}^{2m-2}(s^{2m-u-1}t+1)(s^{2m-u}-1) ~ < ~ \prod_{u=m}^{2m-2}(s^{u+1}-1).
\end{equation} 
As $t \leq s^2$, if we can prove that
\begin{equation}\label{2 bis}
\prod_{u=m}^{2m-2}(s^{2m-u+1}+1)(s^{2m-u}-1) ~ < ~ \prod_{u=m}^{2m-2}(s^{u+1}-1)
\end{equation}  
then (\ref{1 bis}) is proved. Put $v := 2m-u$ in the left side of (\ref{2 bis}) and $v := u-m+2$ in the right side. Then (\ref{2 bis}) is turned into the following:
\begin{equation}\label{3 bis}
\prod_{v=2}^{m}(s^{v-1}+1)(s^v-1) ~ < ~ \prod_{v=2}^{m}(s^{v+m-1}-1).
\end{equation}  
Inequality (\ref{3 bis}) fails to hold for $m \leq 5$ but it holds true for $m = 6$, as one can check. We shall now prove by induction that (\ref{3 bis}) holds for any $m \geq 6$, without taking care of the fact that anyway $m \leq 9$. We must prove that, if 
\begin{equation}\label{4 bis}
\frac{\prod_{v=2}^m(s^{v+1}+1)(s^v-1)}{\prod_{v=2}^m(s^{v+m-1}-1)} ~ < ~ 1,
\end{equation}
then 
\begin{equation}\label{5 bis}
\frac{\prod_{v=2}^{m+1}(s^{v+1}+1)(s^v-1)}{\prod_{v=2}^{m+1}(s^{v+m}-1)} ~ < ~ 1.
\end{equation}
We proceed as follows. 
\[\frac{\prod_{v=2}^{m+1}(s^{v+1}+1)(s^v-1)}{\prod_{v=2}^{m+1}(s^{v+m}-1)} ~ = ~ \frac{\prod{v=2}^{m+1}(s^{v+1}+1)(s^v-1)}{\prod_{v=3}^{m+2}(s^{v+m-1}-1)} ~ =\]
\[= ~  \frac{\prod_{v=2}^{m}(s^{v+1}+1)(s^v-1)}{\prod_{v=2}^{m}(s^{v+m-1}-1)}\cdot\frac{(s^{m+2}+1)(s^{m+1}-1)^2}{(s^{2m}-1)(s^{2m+1}-1)}.\]
The first of these two ratios is less than 1 by the inductive hypothesis (\ref{4 bis}) while the second ratio is obviously less than $1$, since $m \geq 6$. 

So, (\ref{1 bis}) holds, provided that $m \geq 6$. If $m \geq 6$ then $m < i_{\mathrm{max}}-1$ if $t \leq s$ and $m < i_{\mathrm{max}}$ when $t > s$.   Hence $|\Delta_m| > |\Delta_{m-1}|$. Thus, if $m \geq 6$ then $|\Delta_{m-1}| > |\Delta_{n-2}|$, as claimed. 

However $m = n/2 \leq i_{\mathrm{max}}-1$ (recall that $m \geq 6$ by assumption). Therefore $|\Delta_m| > |\Delta_{m-1}|$, by Proposition \ref{maximum}. 
It follows that, if $m \geq 6$, then $|\Delta_m| > |\Delta_{m-1}| > |\Delta_{n-2}|$. This also fixes the case $i = m$, provided that $m \geq 6$.   
 
The cases where $2\leq m \leq 5$ remain to consider. Also, $i = m-1$ or $i = m$, but $m \geq 3$ for $i = m$, since we also want $i < n-2$.
The situation for $i = m-1$ is summarized in the next table, where the rows and the columns correspond to the possible values of $m$ and the possible choices of $t$ respectively and the entries $>$ or $<$ stand for $|\Delta_{m-1}| > |\Delta_{n-2}|$ and $|\Delta_{m-1}| < |\Delta_{n-2}|$, respectively. 
\[(n,i) = (2m, m-1) \hspace{10 mm}
\begin{array}{c|c|c|c|c|c|}
{} & 1 & s^{1/2} & s & s^{3/2} & s^2 \\
\hline
2 & < & < & < & < & < \\
\hline
3 & > & < & < & < & < \\
\hline
4 & > & > & > & < & < \\
\hline
5 & > & > & > & > & < \\
\hline
\end{array}\] 
Proposition \ref{lemma com 4} can be exploited to save a few checks. Indeed, according to it, we must put $<$ in the first row of the above table for $t > 1$ and in the second row for $t > s$. 

Let now $(n,i) = (2m, m)$. According to Proposition \ref{lemma com 3}, if $t < s^2$ then $|\Delta_i| \leq |\Delta_{n-2}|$ only if  $i < n/2$. So, $t = s^2$ is the only case to inspect. It turns out that in this case $|\Delta_{m}| > |\Delta_{n-2}|$ for $m = 3, 4, 5$. 

\medskip

\noindent
{\bf Odd case: $n = 2m+1$.} Now $i = m$. The first inequality of (\ref{com spec 1}) now yields $m  \leq  2\log_s t + 3$. So, $m \leq 7$. As in the even case, one can prove that $|\Delta_m| > |\Delta_{n-2}|$ when $6 \leq m \leq 7$.  

As $m = i < n-2 = 2m-1$, necessarily $m \geq 2$. The cases where $2\leq m \leq 5$ remain to consider. The situation is summarized in the next table. The symbol $=$ occurs in it. It means that, in that particular case,  $|\Delta_{m}| = |\Delta_{n-2}|$. 
\[(n,i) = (2m+1, m) \hspace{10 mm}
\begin{array}{c|c|c|c|c|c|}
{} & 1 & s^{1/2} & s & s^{3/2} & s^2 \\
\hline
2 & > & > & = & < & < \\
\hline
3 & > & > & > & > & < \\
\hline
4 & > & > & > & > & > \\
\hline
5 & > & > & > & > & > \\
\hline
\end{array}\] 
Note that, since $i \leq (n-2+\log_st)/2$ by Proposition \ref{lemma com 3}, we know in advance that $>$ occurs at all entries in the first two columns of the above table. Proposition \ref{lemma com 4} can also be used to save a couple of checks. \hfill $\Box$   

\begin{note}
{\rm The case $(n,i,j) = (5,2,3)$ with $s = t$ is a very special one. It corresponds to case considered in the second part of Proposition \ref{maximum}, where for $n = 5$ the number $|\Delta_2| = |\Delta_3|$ is the maximum of the function $\delta_{5,s,s}$.}
\end{note}
\begin{note}
{\rm In the special case of Proposition \ref{res com 2}, in spite of the fact that $|\Delta_2| = |\Delta_3|$, the graphs $(\Delta_2,\perp)$ and $(\Delta_3,\perp)$, as well as the graphs $(\Delta_2, \perp\cup\approx)$ and $(\Delta_3, \perp\cup\approx)$ and the graphs $(\Delta_2, \perp_{\mathrm{max}})$ and  $(\Delta_3, \perp_{\mathrm{max}})$, have different degrees. Explicitly, let $\kappa_i$, $\chi_i$, $\xi_i$ be the degrees of $(\Delta_i,\perp)$, $(\Delta_i, \perp\cup\approx)$ and $(\Delta_i, \perp_{\mathrm{max}})$ respectively, $i = 2, 3$. In the case under consideration we have}
\[\begin{array}{ccl}
\kappa_2 & = & s(s+1)(s^2+1)(s^3+1)(s^2+s+1), \\
\kappa_3 & = & s(s+1)^2(s^2+1);\\
\chi_2 & = & \kappa_2 + s^5(s^3+s+1),\\
\chi_3 & = & \kappa_3 +  s^3(s^2+1)(s+1);\\
\xi_2 & = & s^4(s+1)(s^2+1)(s^2+s+1),\\
\xi_3 & = & s(s+1)^2(s^2+1) ~ = ~ \kappa_3. 
\end{array}\]
{\rm Clearly, $\kappa_2 > \kappa_3$, $\chi_2 > \chi_3$ and $\xi_2 > \xi_3$. Accordingly, $(\Delta_2,\perp)\not\cong(\Delta_3,\perp)$, $(\Delta_2, \perp\cup\approx)\not\cong(\Delta_3, \perp\cup\approx)$ and $(\Delta_2, \perp_{\mathrm{max}})\not\cong(\Delta_3, \perp_{\mathrm{max}})$, as claimed in Theorem \ref{theo 2}. We shall discuss the degrees $\kappa_i$, $\chi_i$ and $\xi_i$ in a general setting in Section \ref{degree sec}.}
\end{note}

\section{Degrees of $(\Delta_i, \perp)$, $(\Delta_i, \perp\cup\approx)$ and $(\Delta_i, \perp_{\mathrm{max}})$}\label{degree sec}

We denote by $\kappa_i$, $\chi_i$ and $\xi_i$ the degrees of the graphs $(\Delta_i,\perp)$, $(\Delta_i, \perp\cup\approx)$ and $(\Delta_i,\perp_{\mathrm{max}})$, respectively. 

Clearly, $\kappa_i < \chi_i$ for every $i  = 0, 1,..., n-1$. Recall that $(\Delta_{n-1}, \perp)$ and $(\Delta_{n-1},\perp_{\mathrm{max}})$ have no edges ($\kappa_{n-1} = \xi_{n-1}=  0$) while $(\Delta_0, \perp\cup\approx)$ is a complete graph ($\chi_0 = |\Delta_0|-1$). So, $\kappa_i > \kappa_{n-1}$ and $\xi_i > \xi_{n-1}$ for every $i < n-1$ while $\chi_i < \chi_0$ for every $i > 0$. Morever $\xi_i = 0$ if $i < (n-2)/2$. Indeed, if two distinct $i$-spaces $I$ and $I'$ span a maximal singular subspace then $\dm(I\cap I') = 2i-n+1 \geq -1$. Consequently $i \geq (n-2)/2$.  

In this section we shall compute $\kappa_i$ and $\xi_i$ for $i < n-1$ and $\chi_i$ for $i > 0$ and, under suitable hypotheses on $i$ or $j$, we shall compare $\kappa_i$ and $\kappa_j$ as well as $\chi_i$ and $\chi_j$ and $\xi_i$ and $\xi_j$ for $i < j$. 

\subsection{The degree $\kappa_i$ of $(\Delta_i,\perp)$}

\subsubsection{Computing $\kappa_i$}

Given two distinct subspaces $I, I' \in \Delta_i$ with $I\perp I'$, let $k := \dm(I\cap I')$. Then $n-1 \geq \dm(\langle I, I'\rangle) = 2i-k$, namely $k \geq  2i-n+1$. However $k \geq -1$. So, 
\[i-1 ~ \geq ~ k ~ \geq ~ k_{\mathrm{min}}(i) ~ := ~ \mathrm{max}(2i-n+1, -1).\]
Obviously,
\[k_{\mathrm{min}}(i) ~ = ~ \left\{\begin{array}{cl}
-1 & \mbox{if } 2i-n+1 \leq -1, \mbox{ namely } i \leq (n-2)/2; \\
2i-n+1 & \mbox{if } 2i-n+1 \geq -1, \mbox{ namely } i \geq (n-2)/2.
\end{array}\right.\]
Given $I \in \Delta_i$, for every $k = k_{\mathrm{min}}(i), k_{\mathrm{min}}(i)+1, ..., i-1$, the number of $k$-subspaces of $I$ is
\[A_{i,k}~:=~ \prod_{u=0}^k\frac{s^{i+1-u}-1}{s^{k+1-u}-1} ~ = ~ \prod_{v=1}^{k+1}\frac{s^{i-k+v}-1}{s^v-1} ~~~ (= 1 \mbox{ if } k = -1).\]
The number of maximal singular subspaces containing $I$ is
\[B_i ~ := ~ \prod_{u=0}^{n-i-2}(s^ut+1) ~~~ (= 1 \mbox{ if } i = n-1).\]
Given a maximal singular subspace $M \supseteq I$ and a $k$-subspace $K \subsetneq I$, let $C_{i,k}$ be the number of $i$-subspaces $I' \subseteq M$ such that $I\cap I' = K$. A standard computation yields
\[C_{i,k} ~ = ~ s^{(i-k)^2}\cdot\prod_{u=0}^{i-k-1}\frac{s^{n-i-1-u}-1}{s^{i-k-u}-1} ~ = ~ s^{(i-k)^2}\cdot\prod_{v=1}^{i-k}\frac{s^{n-2i+k-1+v}-1}{s^{v}-1}.\]
If $I'$ is such an $i$-subspace, then $\dm(\langle I, I\rangle) = 2i-k$, hence $\langle I, I'\rangle$ is contained in exactly
\[D_{i,k} ~ := ~ \prod_{u=0}^{n-2i+k-2}(s^ut + 1) ~~~ (= 1 \mbox{ if } 2i-k = n-1)\] 
maximal singular subspaces. Clearly, the number of $i$-subspaces $I'\perp I$ such that $\dm(I\cap I') = k$ is
\[\kappa_{i,k} ~ = ~ \frac{A_{i,k}B_iC_{i,k}}{D_{i,k}} ~ = ~ F_{i,k}G_{i.k}\]
where 
\[F_{i,k} ~ := ~ \frac{B_i}{D_{i,k}} ~ = ~ \frac{\prod_{u=0}^{n-i-2}(s^ut+1)}{\prod_{u=0}^{n-2i+k-2}(s^ut+1)} ~ = ~ \prod_{u=n-2i+k-1}^{n-i-2}(s^ut+1),\]
\[G_{i,k} ~ := ~ A_{i,k}C_{i,k} ~ = ~ s^{(i-k)^2}\cdot\prod_{v=1}^{k+1}\frac{s^{i-k+v}-1}{s^v-1}\cdot\prod_{v=1}^{i-k}\frac{s^{n-2i+k-1+v}-1}{s^v-1}.\]
Finally,
\begin{equation}\label{kappa 1}
\kappa_i ~ = ~ \sum_{k = k_{\mathrm{min}}(i)}^{i-1}\kappa_{i,k} ~ = ~ \sum_{k = k_{\mathrm{min}}(i)}^{i-1}F_{i,k}G_{i,k}.
\end{equation}  

\subsubsection{Comparing $\kappa_i$ with $\kappa_j$} 

We firstly consider the case $(n-2)/2 \leq i$. 

\begin{lemma}\label{kappa com 1}
If $(n-2)/2\leq i < n-1$ then $\kappa_i > \kappa_{i+1}$.
\end{lemma}
{\bf Proof.} We have $k_{\mathrm{min}}(i) = 2i-n+1 \geq -1$ and $k_{\mathrm{min}}(i+1) = 2(i+1)-n+1 = 2i-n+3 \geq 1$. According to (\ref{kappa 1}) we have
\[\kappa_i - \kappa_{i+1} ~ = ~ F_{i,2i-n+1}G_{i,2i-n+1} + F_{i,2i-n+2}G_{2i-n+2} - F_{i+1,i}G_{i+1,i} + \]
\[\hspace{20 mm} + \sum_{k = 2i-n+3}^{i-1}(F_{i,k}G_{i,k}-F_{i+1,k}G_{i+1,k}).\]
If we prove the following, we are done: 
\begin{equation}\label{kappa com 1}
F_{i,2i-n+1}G_{i,2i-n+1} + F_{i,2i-n+2}G_{2i-n+2} ~ > ~  F_{i+1,i}G_{i+1,i},
\end{equation}
\begin{equation}\label{kappa com 2}
F_{i,k}G_{i,k} ~ > ~ F_{i+1,k}G_{i+1,k}, \hspace{5 mm}\mbox{for }~ 2i-n+3 \leq k \leq i-1. 
\end{equation}
Inequality (\ref{kappa com 1}) is straightforward. We leave details for the reader (recall that $(n-2)/2 \leq i \leq n-2$, by assumption). In order to prove (\ref{kappa com 2}) we compute the ratio $F_{i,k}G_{i,k}/F_{i+1,k}G_{i+1,k}$. By routine computations we obtain that
\[\frac{F_{i,k}G_{i,k}}{F_{i+1,k}G_{i+1,k}} ~ = ~ \frac{s^{n-2i-1}t+1}{s^{n-2i+k-2}t+1}\cdot\frac{s^{n-2i-2}t+1}{s^{n-2i+k-3}t+1}\cdot s^{2i(k-1)-1}\cdot\]
\[\hspace{27 mm}\cdot\frac{s^{i-k+1}-1}{s^{i+2}-1}\cdot\frac{s^{n-i-1}-1}{s^{n-2i+k-1}-1}\cdot\frac{s^{n-i-2}-1}{s^{n-2i+k-2}-1}.\]
As $k < i$, we have
\[\frac{s^{n-i-1}-1}{s^{n-2i+k-1}-1} ~ > ~ \frac{s^{n-i-1}}{s^{n-2i+k-1}} ~ = ~ s^{i-k},\]
\[\frac{s^{n-i-2}-1}{s^{n-2i+k-2}-1} ~ > ~ \frac{s^{n-i-2}}{s^{n-2i+k-2}} ~ = ~ s^{i-k}.\]
Moreover, since $i \geq (n-2)/2$ and $k \geq 2i-n+3$, we also have $k \geq 1$. Hence
\[\frac{s^{n-2i-1}t+1}{s^{n-2i+k-2}t+1} ~ \geq ~ \frac{s^{n-2i-1}t}{s^{n-2i+k-2}t} ~ = ~ s^{1-k},\]
\[\frac{s^{n-2i-2}t+1}{s^{n-2i+k-3}+1} ~ \geq ~ \frac{s^{n-2i-2}t}{s^{n-2i+k-3}t} ~ = ~ s^{1-k}.\]
Therefore
\[\frac{F_{i,k}G_{i,k}}{F_{i+1,k}G_{i+1,k}} ~ > ~ s^{2-2k}\cdot s^{2i(k-1)-1}\cdot\frac{s^{i-k+1}-1}{s^{i+2}-1}\cdot s^{2i-2k},\] 
namely:
\begin{equation}\label{kappa com 3}
\frac{F_{i,k}G_{i,k}}{F_{i+1,k}G_{i+1,k}} ~ > ~ s^{2ik-4k+1}\cdot\frac{s^{i-k+1}-1}{s^{i+2}-1}.
\end{equation} 
It is easy to check that the right side of (\ref{kappa com 3}) is greater than $1$. Inequality (\ref{kappa com 2}) is proved. Thus, the lemma is proved as well. \hfill $\Box$\\

By Lemma \ref{kappa com 1} we immediately obtain the following:

\begin{prop}\label{kappa com 4}
Let $(n-2)/2\leq i < j \leq n-1$. Then $\kappa_i > \kappa_j$.
\end{prop}

In view of Proposition \ref{kappa com 4}, if $m \in \{0, 1,..., n-1\}$ is such that $\kappa_\nu$ is maximal, then $m \leq (n-2)/2$. Regretfully, we are presently unable to determine $m$, but we can prove that $m \geq (n-5+\log_st)/3$.

\begin{lemma}\label{kappa com 5}
Suppose $n \geq 5$ and let $i < (n-5+\log_st)/3$. Then $\kappa_i < \kappa_{i+1}$.
\end{lemma}
{\bf Proof.} Now $k_{\mathrm{min}}(i) = k_{\mathrm{min}}(i+1) = -1$. According to (\ref{kappa 1}) we have
\[\kappa_{i+1}- \kappa_{i} ~ = ~ F_{i+1,i}G_{i+1,i} + \sum_{k = -1}^{i-1}(F_{i+1,k}G_{i+1,k}-F_{i,k}G_{i,k}).\]
If we prove that
\begin{equation}\label{kappa com 6}
F_{i+1,k}G_{i+1,k} ~ \geq ~ F_{i,k}G_{i,k} \hspace{5 mm}\mbox{for }~ -1 \leq k \leq i-1 
\end{equation}
then we are done. By routine computations we get
\[\frac{F_{i+1,k}G_{i+1,k}}{F_{i,k}G_{i,k}} ~ = ~ s^{2i-2k}(s^{n-2i+k-3}t+1)\cdot\frac{s^{n-2i+k-2}t+1}{s^{n-i-2}t+1}\cdot\frac{s^{i+2}-1}{s^{i-k+1}-1}\cdot\]
\[\hspace{27 mm} \cdot \frac{(s^{n-2i+k-1}-1)(s^{n-2i+k-2}-1)}{(s^{n-i-1}-1)(s^{n-2-i}-1)} ~ = \]
\[= ~ s^{2i-2k}(s^{n-2i+k-3}t+1)\cdot\frac{s^{n-2i+k-2}t+1}{s^{n-i-2}t+1}\cdot\frac{s^{i+2}-1}{s^{i-k+1}-1}\cdot\]
\[\hspace{5 mm} \cdot \frac{s^{2n-4i+2k-3}-s^{n-2i+k-1}-s^{n-2i+k-2}+1}{s^{2n-2i-3}-s^{n-i-1}-s^{n-i-2}+1}.\]
Clearly, $s^{n-2i+k-2} < s^{n-i-2}$, $s^{i+2}\geq s^{i-k+1}$. and
\[0  <  s^{2n-4i+2k-3}-s^{n-2i+k-1}-s^{n-2i+k-2} ~ < ~ s^{2n-2i-3}-s^{n-i-1}-s^{n-i-2}\] 
since $i < (n-3)/3$. Accordingly,
\[\frac{s^{n-2i+k-2}t+1}{s^{n-i-2}t+1} ~ > ~ \frac{s^{n-2i+k-2}t}{s^{n-i-2}t} ~ = ~ s^{k-i},\]
\[\frac{s^{i+2}-1}{s^{i-k+1}-1} ~ \geq ~ \frac{s^{i+2}}{s^{i-k+1}} ~ = ~ s^{k+1},\]
\[\frac{s^{2n-4i+2k-3}-s^{n-2i+k-1}-s^{n-2i+k-2}+1}{s^{2n-2i-3}-s^{n-i-1}-s^{n-i-2}+1} ~ >\]
\[> ~ \frac{s^{2n-4i+2k-3}-s^{n-2i+k-1}-s^{n-2i+k-2}}{s^{2n-2i-3}-s^{n-i-1}-s^{n-i-2}} ~ =\]
\[= ~ \frac{s^{n-3i+k-1}-s^{k-i+1}-s^{k-i}}{s^{n-i-1}-s-1}.\]
Consequently,
\[\frac{F_{i+1,k}G_{i+1,k}}{F_{i,k}G_{i,k}} ~ > ~ (s^{n-2i+k-3}t+1)\cdot\frac{s^{n-2i+2k}-s^{k+2}-s^{k+1}}{s^{n-i-1}-s-1}.\]
Exploiting the assumptions $n \geq 5$ and $i < (n-5+\log_st)/3$ it is straightforward to check that 
\[(s^{n-2i+k-3}t+1)\cdot\frac{s^{n-2i+2k}-s^{k+2}-s^{k+1}}{s^{n-i-1}-s-1} ~ > ~ 1.\]
for any $k$. So, (\ref{kappa com 6}) is proved. \hfill $\Box$

\begin{note}
{\rm When $n = 5$ or $n = 6$ with $t \in \{1, s^{1/2}, s, s^{3/2}\}$ or $n = 7$ with $t \in \{1, s^{1/2}\}$ the hypothesis that $i < (n-5+\log_st)/3$ forces $i = 0$. In this cases Lemma \ref{kappa com 5} only says that $\kappa_0 < \kappa_1$. The same when $n = 4$ with $t\in \{s, s^{3/2}, s^2\}$. When $n = 4$ with $t \in \{1, s^{1/2}\}$ or $n = 3$ then $(n-5+\log_st)/3 < 0$, whence the hypothesis $i < (n-5+\log_st)/3$ is vacuous. The next lemma covers these extremal cases.}
\end{note}

\begin{lemma}\label{kappa com 7}
If $n = 4$ then $\kappa_0 < \kappa_1 > \kappa_2 > \kappa_3$. If $n = 3$ then $\kappa_0 > \kappa_1 > \kappa_2$.
\end{lemma}
{\bf Proof.} Straightforward verification. Note that the inequalities $\kappa_1 > \kappa_2 > \kappa_3$ (for $n = 4$) and $\kappa_1 > \kappa_2$ (for $n = 3$)
also follow from Lemma \ref{kappa com 1} noticing that $(n-2)/2 \leq 1$ when $n \leq 4$. \hfill $\Box$ \\

By Lemmas \ref{kappa com 5} and \ref{kappa com 7} we obtain the following:

\begin{prop}\label{kappa com 8}
Suppose that $n \geq 5-\log_st$ and let $i < j \leq  \lceil(n-5+\log_st)/3\rceil$. Then $\kappa_i < \kappa_j$.
\end{prop}

\subsection{The degree $\chi_i$ of $(\Delta_i, \perp\cup\approx)$} 

Given $I \in \Delta_i$, let $\lambda_i$ be the number of $i$-subspaces $I'\in \Delta_i$ such that $\dm(I\cap I') = i-1$ but $I\not\perp I'$. Given an $(i-1)$-subspace $K\subsetneq I$, the subspace $I$ is one of the $(s^{n-i-1}t+1)(s^{n-i}-1)/(s-1)$ points of the polar space $\mathrm{Res}_\Delta(K)$. This polar space contains exactly $s^{2n-2i-2}t$ points non-collinear with $I$. Since $I$ contains exactly $(s^i-1)/(s-1)$ subspaces of dimension $i-1$, it turns out that
\[\lambda_i ~ = ~ s^{2n-2i-2}t\cdot\frac{s^{i+1}-1}{s-1}.\]
Clearly, $\chi_i = \kappa_i + \lambda_i$. So:
\begin{equation}\label{chi 1}
\chi_i ~ = ~ \kappa_i + s^{2n-2i-2}t\cdot\frac{s^{i+1}-1}{s-1}.
\end{equation} 
It is easily seen that $\lambda_i > \lambda_{i+1}$ for every $i < n-1$. So, by Proposition \ref{kappa com 4} and (\ref{chi 1}) we obtain the following

\begin{prop}\label{chi 2}
Let $(n-2)/2\leq i < j \leq n-1$. Then $\chi_i > \chi_j$.
\end{prop}

\subsection{The degree $\xi_i$ of $(\Delta_i, \perp_{\mathrm{max}})$} 

As previously noticed, $\xi_i = 0$ if either $i < (n-2)/2$ or $i = n-1$. So, assume that $(n-2)/2 \leq i < n-1$. 

Let $I\in \Delta_i$. For every $M\in \Delta_{n-1}$ containing $I$, the $i$-subspaces  $I' \neq I$ of $M$ such that $I'\perp_{\mathrm{max}} I$ bijectively correspond to the pairs $(K, H)$ where $K$ is a $(2i-n+1)$-subspace of $I$ and $H$ is a complement of the quotient space $I/K$ in $M/K$. Accordingly, $\xi_i = A_iB_iC_i$ where $A_i$ is the number of maximal singular subspaces containing $I$, $B_i$ is the number of $(2i-n+1)$-subspaces of $I$ and $C_i$ is the number of complements of $I/K$ in $M/K$, for a given $M\in \Delta_{n-1}$ containing $I$ and a given $(2i-n+1)$-subspace $K$ of $I$. Explicitly:
\[A_i ~ = ~ \prod_{u=0}^{n-i-2}(s^ut+1),\]
\[B_i ~ = ~ \frac{\prod_{u=0}^{2i-n+1}(s^{i+1-u}-1)}{\prod_{u=0}^{2i-n+1}(s^{2i-n+2-u}-1)} ~ = ~  \frac{\prod_{u=n-i}^{i+1}(s^u-1)}{\prod_{u=1}^{2i-n+2}(s^u-1)} \hspace{2 mm} (= 1 \mbox{ if } i = \frac{n-2}{2}), \]
\[C_i ~ = ~ s^{(n-i-1)^2}.\]
Clearly, $\xi_i = A_iB_iC_i$. Explicitly: 
\begin{equation}\label{xi 1}
\xi_i ~ = ~  s^{(n-i-1)^2}\cdot\prod_{u=0}^{n-i-2}(s^ut+1)\cdot\frac{\prod_{u=n-i}^{i+1}(s^u-1)}{\prod_{u=1}^{2i-n+2}(s^u-1)}.
\end{equation}

\begin{prop}\label{xi 2}
Let $(n-2)/2\leq i < j \leq n-1$. Then $\xi_i > \xi_j$.
\end{prop}
{\bf Proof.} As $\xi_{n-1} = 0$ we can assume that $j < n-1$. Accordingly, $i < n-2$. For $(n-2)/2 \leq i < n-2$, (\ref{xi 1}) yields
\[\frac{\xi(i+1)}{\xi(i)} ~ = ~ \frac{(s^{n-i-1}-1)(s^{i+2}-1)}{s^{2n-2i-3}(s^{n-i-2}t+1)(s^{2i-n+4}-1)(s^{2i-n+3}-1)}.\]
By standard manipulations one can check that 
\[(s^{n-i-1}-1)(s^{i+2}-1)  <  s^{2n-2i-3}(s^{n-i-2}t+1)(s^{2i-n+4}-1)(s^{2i-n+3}-1)\]
for any $i \in \{\lceil (n-2)/2\rceil,..., n-3\}$. The conclusion follows. \hfill $\Box$ 

\subsection{Proof of Theorem \ref{theo 2}}\label{proof}

By Propositions \ref{kappa com 4}, \ref{chi 2} and \ref{xi 2}, if $(n-2)/2\leq i < j \leq n-1$ then $\kappa_i > \kappa_j$, $\chi_i > \chi_j$ and $\xi_i > \xi_j$. The conclusions of Theorem \ref{theo 2} follow.  

\begin{conj}
For $0 \leq i < n$, let $\theta_i$ be any of $\kappa_i, \chi_i$ or $\xi_i$. If $0 \leq i < j < n$, then $(|\Delta_i|, \theta_i)\neq (|\Delta_j|, \theta_j)$.
\end{conj}

\section{The degrees of $(\Delta_i, \approx)$ and $(\Delta_i, \approx\cap\perp)$}

For the sake of completeness, we also give the degrees of $(\Delta_i, \approx)$ and $(\Delta_i, \approx\cap\perp)$. Let $\mu_i$ and $\nu_i$ be the degrees of  $(\Delta_i, \approx)$ and $(\Delta_i, \approx\cap\perp)$ respectively.

Let $I \in \Delta_i$. The subspace $I$ contains exactly $A_i = (s^{i+1}-1)/(s-1)$ hyperplanes and, if $H$ is one of them, its residue $\mathrm{Res}_\Delta(H)$ contains exactly $B_i = (s^{n-i-1}t+1)(s^{n-i}-1)/(s-1)$ points. The subspace $I$ is one of them. Therefore 
\[\mu_i ~  = ~  A_i(B_i-1) ~ = ~  s(s^{2n-2i-2}t - s^{n-i-2}t + s^{n-i-1} - 1)\cdot\frac{s^{i+1}-1}{(s-1)^2}.\]
The polar space $\mathrm{Res}_\Delta(H)$ contains exactly 
\[C_i ~ = ~ (s^{n-i-2}t+1)\cdot\frac{(s^{n-i-1}-1}{s-1}\cdot s\]
points collinear with $I$ (and different from $I$). Therefore 
\[\nu_i ~ = ~ A_iC_i ~ = ~  s(s^{n-i-2}t+1)\cdot\frac{(s^{n-i-1}-1)(s^{i+1}-1)}{(s-1)^2}.\]
It is easily checked that if $i < n-1$ then $\mu_i > \mu_{i+1}$ and $\nu_i > \nu_{i+1}$. So, if $i < j$ then  $\mu_i > \mu_j$ and $\nu_i > \nu_j$.

\end{document}